\documentclass[graybox]{svmult}

\usepackage{type1cm}
\usepackage[numbers]{natbib}
\usepackage{makeidx}
\usepackage{graphicx}
\usepackage{multicol}
\usepackage[bottom]{footmisc}
\usepackage{graphicx}
\usepackage{mathtools}
\usepackage{hyperref}
\usepackage{pgf,tikz}
\usepackage{mathrsfs}
\usepackage{float}
\usepackage{csvsimple}
\usepackage{enumerate}
\usepackage{multirow}
\usepackage{siunitx}
\usepackage{pgfplots}
\usepackage{pgfplotstable}
\usepackage[noend]{algpseudocode}
\usepackage{graphicx,bm,xcolor}
\usepackage{algorithm}
\usepackage{pdfpages}
\usepackage{stmaryrd}
\usepackage{siunitx}
\usepackage{hyperref}
\usetikzlibrary{arrows}
\usetikzlibrary{math} 
\usepackage{caption}
\usepackage{t1enc}
\usepackage[polish,german,english]{babel}
\usepackage{verbatim} 
\usepackage{xcolor}

\usepackage{url}

\newif\ifMAKEPICS
\MAKEPICStrue

\ifMAKEPICS
\usepackage[cleanup,subfolder]{gnuplottex}
\usepackage{xparse}

\ExplSyntaxOn
\DeclareExpandableDocumentCommand{\convertlen}{ O{cm} m }
{
	\dim_to_decimal_in_unit:nn { #2 } { 1 #1 } cm
}
\ExplSyntaxOff
\fi

\usepackage{newtxtext}
\usepackage[varvw]{newtxmath}

\newcommand*{\totaldt}{\frac{\mathrm{d}}{\mathrm{d}t}}

\newcommand*{\tsep}{t_{\mathrm{separate}}}
\newcommand*{\tdur}{t_{\mathrm{dur}}}
\newcommand*{\tcont}{t_{\mathrm{contact}}}

\begin{document}

\title*{Goal Oriented Adaptive Space Time Finite Element Methods Applied to Touching Domains}
\titlerunning{G.O.A. Space-time FEM Applied to Touching Domains}
\author{Bernhard Endtmayer\orcidID{0000-0003-0647-7371} and\\ Andreas Schafelner\orcidID{0000-0001-7997-1165}}
\institute{B. Endtmayer \at Leibniz Universit\"at Hannover, Institut f\"ur Angewandte Mathematik, Welfengarten 1, 30167 Hannover, Germany
\at Cluster of Excellence PhoenixD (Photonics, Optics, and Engineering -- Innovation Across Disciplines), Leibniz Universit\"at Hannover, Germany,  \email{endtmayer@ifam.uni-hannover.de}
\and A. Schafelner \at Institute for Numerical Mathematics, Johannes Kepler University Linz, Altenberger Stra\ss{}e 69, 4040 Linz, Austria \email{andreas.schafelner@jku.at}}
\maketitle

\abstract{We consider goal-oriented adaptive space-time finite-element discretizations of the parabolic heat equation on completely unstructured simplicial space-time meshes.
In some applications, we are interested
in an accurate computation of some possibly nonlinear functionals at the solution, so called goal functionals.
This motivates the use of adaptive mesh refinements driven by the dual-weighted residual (DWR) method.
The DWR method requires the numerical solution of a linear adjoint problem 
that provides the sensitivities for the mesh refinement. This can be done 
by means of the same full space-time finite element discretization as used 
for the primal linear problem.
The numerical experiment presented demonstrates that this goal-oriented, full space-time  
finite element solver efficiently provides accurate numerical results 
for a model problem with moving domains and a linear goal functional, where we know the exact value.}

\section{Introduction}\label{sec:intro}

We consider the following parabolic evolution equation on moving domains: find $u$ such that
\begin{gather}\label{eq:model problem}
    \frac{\mathrm{d}}{\mathrm{d}t} u -\Delta_x u = 0\ \text{in } Q, \\
    \nabla_x u \cdot \mathbf{n} = 0\ \text{on } \Sigma,\text{ and }u = u_0\ \text{ on } \Sigma_0,\label{eq:model problem:ibc}
\end{gather}
with the material derivative $ \frac{\mathrm{d}}{\mathrm{d}t} = \partial_t + (v \cdot\nabla_x) $, where $ v = \partial_t \Phi $ and $ \Phi$ describes the movement of the the domain(s), the spatial Laplacian $ \Delta_x$, and $ Q = \{ (x,t) \vcentcolon t\in(0,T),\ x\in \Omega(t) \} $, $ \Sigma = \{ (x,t) \vcentcolon t\in(0,T),\ x\in\partial\Omega(t) \} $, $ \Sigma_0 = \{(x,0) \vcentcolon x\in\Omega(0) \}$, where $ \Omega(t) = \{ \Phi(x,t) \vcentcolon x\in\Omega\} $ for $t\in(0,T)$.

We are however not interested in the solution $u$ of \eqref{eq:model problem}--\eqref{eq:model problem:ibc}, but instead in the value of a certain Quantity of Interest (QoI), which is typically represented by a possibly nonlinear functional $J(\cdot)$. The QoI is usually some localized object, e.g.\ the mean of $u$ over a certain subregion of $\Omega$ at a certain time $t$. This warrants the use of adaptive mesh refinement to reduce the overall computational cost when aiming for a sufficiently high local resolution of the QoI. The standard approach for the solution of this class of problems is the method of lines, where we perform semi-discretization in space and time separately, 
where the {goal oriented} adaptivity introduces additional challenges, e.g.\ we have to store the solution at all discretization points in time; see e.g.\ \cite{besier2010goal,failer2018adaptive,munoz2019explicit,fischer2023adaptive} and the references therein. One possible way to overcome this requires checkpoint techniques as presented in \cite{ThesisMeisrimel}.
In this paper, we instead treat time as just another variable and perform an all-at-once discretization of the complete space-time cylinder $Q$. While this increases the dimension of our problem by one, it also enables us to simultaneously refine in space and time, and we have by design access to the solution over the whole time interval. Moreover, since we know the movement of the spatial domain in advance, the space-time domain is fixed; see e.g.\ Fig.~\ref{fig:space-time domain}.
For more details on goal oriented space-time adaptivity for parabolic evolution problems, see \cite{endtmayer2023goaloriented} and the references therein.

\tikzmath{\tcontact=.45; \tduration=.4; \tseparate=\tcontact+\tduration; \tend=\tcontact+\tduration+\tcontact;\x1=2;}
\begin{figure}
    \centering%
    \begin{tikzpicture}[scale=1.5]
        \draw[gray,<->] (-1.2*\x1,0) -- (0,0) -- (1.2*\x1,0) node[anchor=north west,scale=.75] {$x$};
        \draw[gray,->] (0,0) -- (0,1.2*\tend) node[anchor=south west,scale=.75] {$ t $};
        \draw (-2,0) -- node[anchor=north,scale=.75] {$ \Omega_{1}$ } (-1,0) -- ++(1,\tcontact) -- (1,0) -- node[anchor=north,scale=.75] {$ \Omega_{2}$ } (2,0) -- ++(-1,\tcontact) -- ++(0, \tduration) -- (2,\tend) -- ++(-1,0) -- (0, \tseparate) -- (-1, \tend) -- (-2,\tend) -- (-1, \tseparate) -- (-1,\tcontact) -- cycle;
        \draw[gray,densely dashed] (-2,\tcontact) node[anchor=north east,scale=.75] {$ \tcont $} -- ++(4,0);
        \draw[gray,densely dashed] (-2,\tseparate) node[anchor=south east,scale=.75] {$ \tsep $} -- ++(4,0);
        \draw[gray,<->] (2,\tcontact) -- node[anchor=west,scale=.75] {$\tdur$}  ++(0,\tduration);
    \end{tikzpicture}
    \caption{Two dimensional space-time domain, with the spatial dimension on the $x$-axis, and the temporal dimension on the $y$-axis.}\label{fig:space-time domain}
\end{figure}

The rest of the paper is structured as follows: in Section~\ref{sec:model problem}, we formally introduce the model problem and the space-time finite element discretization. In Section~\ref{sec:go estimates} we recall basic properties for the goal oriented error estimation technique and apply it to the space-time problem. Finally, in Section~\ref{sec:numerical experiments}, we present a numerical experiment in one space dimension.

\section{The model problem and discretization}\label{sec:model problem}

In this section, we briefly describe the model problem and the corresponding space-time finite element discretization.

\subsection{The model problem}

We now consider \eqref{eq:model problem}--\eqref{eq:model problem:ibc} for the configuration sketched in Fig.~\ref{fig:space-time domain}, i.e.\ we consider two separate domains $ \Omega_1 $ and $\Omega_2$, which move towards each other with constant velocity $v$. Upon contact at $ t = \tcont > 0$, the movement stops and the domains stay in contact for a fixed period of time $ \tdur > 0$ , after which they separate at time $ t = \tsep = \tcont + \tdur$ and move away from each other, again with constant velocity $v$. In this setting, the total time derivative can be written as $ \totaldt = \partial_t + v\,\partial_x $.

By standard arguments we obtain the following weak formulation from \eqref{eq:model problem}--\eqref{eq:model problem:ibc}:
Find $u \in Eu_0+U:=\{Eu_0+\psi: \psi \in U\}$ such that
\begin{equation}
\langle \totaldt u, v \rangle + \langle \nabla_x u, \nabla_x v\rangle = \langle f,v \rangle \qquad \quad \forall v \in V,
\end{equation}
where $U$ and $V$ are reflexive Banach-spaces. 
Here, $ \langle \cdot,\cdot\rangle $ denotes the $V^*\times V $ duality product, and $Eu_0$ solves the problem: Find $Eu_0 \in H^1(Q)$ such that 
\begin{align*}
-\Delta Eu_0 =0 \text{ in } Q, \qquad \text{and}\quad Eu_0=u_0  \text{ on }\Sigma_0, \quad \text{and}\qquad \nabla Eu_0\cdot n=0  \text{ on }\Sigma.
\end{align*}
\subsection{Discretization}

Let us assume $Q$ is polytopic. We decompose the domain $Q$ into non-overlapping shape-regular simplicial finite elements $\Delta \in \mathcal{T}_h$. Then, we define the corresponding finite space 
\begin{equation}
\label{disc:FE-Space}
U_h=V_h=\{v_h\in S^k_h(\overline{Q}): v_h=0\text{ on }\Sigma_0  \},
\end{equation}
where $S^k_h(\overline{Q}):=\{v_h \in \mathcal{C}(\overline{Q}): v_h \in \mathbb{P}^k(\Delta) \}$ is nothing else but the standard finite element space with polynomial degree $k=1,2,3, \ldots$. Here, $\mathbb{P}^k(\Delta)$ describes the polynomials with degree $k$ on the element $\Delta$. This leads to a conforming discretization.
With this we define the discretized problem as: Find $u_h \in E_hu_0+U_h$ such that 
\begin{equation}
\label{disc:Primal-Problem}
\langle \totaldt u_h, v_h \rangle + \langle \nabla_x u_h, \nabla_x v_h\rangle = 0 \qquad \quad \forall v_h \in V_h,
\end{equation}
{where $E_h$ is an extension operator into $S_h^k(Q)$, which evaluates $u_0$ at the essential space-time boundary $ \Sigma_0$ and is zero everywhere else.}
Further information about the {space-time} discretization can be found in \cite{steinbach2015}. 

\section{Goal oriented error estimates}\label{sec:go estimates}
In many applications the solution itself is not of primary interest but one or several quantities of interest $J_i \vcentcolon Eu_0+U \to \mathbb{R}$ which are evaluated at the solution. In this work, we just consider one quantity of interest $J$. For multiple quantities of interest we refer to \cite{HaHou03,Ha08,EnWi17,EnLaNeiWoWi2020,EnLaWi18}. Here, our error estimator should estimate the error in our functional, i. e. $\eta_h \approx J(u)-J(u_h)$. For this we use the DWR (dual weighted residual) method \cite{BeRa01,BaRa03}.
In general, the DWR method requires the solution of the adjoint problem. However in \cite{endtmayer2023goaloriented}, it was shown that just the solution of an discretized adjoint problem is required. The discretized adjoint is given by: Find $z_h \in W_h$ such that:
\begin{equation}
\langle \totaldt v_h, z_h \rangle + \langle \nabla_x v_h, \nabla_x z_h\rangle = \langle J'(u_h), v_h\rangle, \qquad \forall v_h \in W_h,
\end{equation}
where $u_h$ solves \eqref{disc:Primal-Problem} and $J'$ describes the Fr\'echet derivative of the goal functional $J$. For linear quantities of interest the derivative $J'$ does not depend on the solution of the primal problem $u_h$. The space $W_h$ can be the finite element space $V_h$ introduced above or an enriched finite element space $V_h^{(2)}$ with the property $V_h \subset V_h^{(2)} \subset V$. 
In order to construct the error estimator as in \cite{endtmayer2023goaloriented}, the solution of the enriched primal problem is necessary. The enriched primal problem is given by: Find $u_h^{(2)} \in V_h^{(2)}$ such that
\begin{equation}
\langle \totaldt u_h^{(2)}, v_h^{(2)} \rangle + \langle \nabla_x u_h^{(2)}, \nabla_x v_h^{(2)}\rangle = \langle f,v_h^{(2)} \rangle \qquad \quad \forall v_h^{(2)} \in V_h^{(2)}.
\end{equation}
Finally the resulting error estimator $\eta_h$ is given by
\begin{equation}
\eta_h:= \frac{1}{2}\rho(u_h)(z_h^{(2)}-z_h)+\frac{1}{2}\rho^*(u_h,z_h)(u_h^{(2)}-u_h)+ \rho(u_h)(z_h)+\mathcal{R}_h^{(3)(2)}.
\end{equation}
This error estimator proofs to be efficient and reliable as shown in \cite{EnLaWi20,endtmayer2021reliability}. 
In particular for our model problem the resulting terms look as follows:
\begin{align*}
\rho(u_h)(e_h^*)&= \langle f,e_h^*\rangle-\langle\totaldt u_h, e_h^* \rangle - \langle \nabla_x u_h, \nabla_x e_h^*\rangle,\\
\rho^*(u_h,z_h)(e_h)&= \langle J'(u_h),e_h\rangle-\langle \totaldt e_h, z_h \rangle - \langle \nabla_x z_h, \nabla_x e_h\rangle,\\
\rho(u_h)(z_h)&= \langle f,z_h\rangle-\langle\totaldt u_h, z_h \rangle - \langle \nabla_x u_h, \nabla_xz_h\rangle,\\
\mathcal{R}_h^{(3)(2)}&=\frac{1}{2} \int_{0}^{1}s(s-1)J'''(s u_h^{(2)}+(1-s)u_h)(e_h,e_h,e_h) \mathrm{d}s,
\end{align*}
where $J'''$ describes the 3rd Fr\'echet derivative, $e_h:=u_h^{(2)}-u_h$, and $e_h^*:=z_h^{(2)}-z_h$. Of course for linear goal functionals $\mathcal{R}_h^{(3)(2)}=0$ and $\rho(u_h)(z_h^{(2)}-z_h) = \rho^*(u_h,z_h)(u_h^{(2)}-u_h)$.
The localization is based on the partition of unity technique in \cite{RiWi15_dwr}. For further details and information we refer to \cite{endtmayer2023goaloriented}.
The resulting adaptive space-time algorithm is summarized in Algorithm~\ref{alg:Full Adaptive Algorithm}.

\begin{algorithm}
    \caption{The adaptive space-time algorithm \label{alg:Full Adaptive Algorithm}}
    \begin{algorithmic}[1] 
    \Repeat
    \State solve the primal and adjoint problem
    \State solve the enriched primal and adjoint problem 
    \State compute the elementwise contributions via PU-technique
    \State select a set of marked elements $M$ using D\"orfler marking \cite{Doerfler:1996a}, 
    \State $\mathcal{T}_{k+1} \gets \Call{Refine}{\mathcal{T}_k, M} $,
    \State $ k \gets k + 1$,
    \Until some stopping criterion is fulfilled. 
    \end{algorithmic}
\end{algorithm}

\section{Numerical experiments}\label{sec:numerical experiments}

Now let {$Q = \mathring{\overline{M(\Omega_1 \cup \Omega_2)}}  $} with $ \Omega_1 = (-2,-1) $ and $ \Omega_2 = (1,2) $,
where
 \begin{align*}
 M(\Omega)&:=\{m(x,t): x \in \Omega, t \in (0,T)\}, \\
m(x,t)& = \begin{cases}
(x-\frac{20xt}{9|x|},t) &\ t < 0.45,\\
(x+\frac{20xt}{9|x|}-\frac{20x}{9|x|},t) &\ t > 0.55,\\
(x-\frac{x}{|x|},t) & \text{else}.
\end{cases}
\end{align*}
We visualized the domain $Q$ in Figure~\ref{fig:space-time domain}.
 We consider the following initial boundary value problem:
find $u$ such that
\begin{gather*}
    \partial_t u (x,t) + v(x,t)\,\partial_x u(x,t) - \partial_{x^2}^2 u(x,t) = 0, (x,t) \in Q\\
    \partial_x u(x,t)\,n_x = 0, x\in \Sigma \\
    u(x,0) = u_0(x), x \in {\Omega_1}\cup{\Omega_2}
\end{gather*}
where
\begin{equation*}
    v(x,t) = \begin{cases}
        \frac{20}{9} & x < 0\ \text{and}\ t < 0.45,\\
        -\frac{20}{9} & x > 0\ \text{and}\ t < 0.45,\\
        -\frac{20}{9} & x < 0\ \text{and}\ t > 0.55,\\
        \frac{20}{9} & x > 0\ \text{and}\ t > 0.55,\\
        0 & \text{else},
    \end{cases}
    \quad\text{and}\quad
    u_0(x) = \begin{cases}
        2 & x\in\Omega_1, \\
        0 & \text{else}.
    \end{cases}
\end{equation*}

In order to test the adaptive Algorithm~\ref{alg:Full Adaptive Algorithm}, we are interested in the following Quantity of Interest
\begin{align*}
    J(u) &\coloneqq \int_{1}^{2}\! u(x,T)\;\mathrm{d}x = g(\tdur) \approx \num[round-precision=15,round-mode=places]{0.3568234004524540428094709313813},\\
    g(\tdur)&\coloneqq\frac{1}{\sqrt{\pi \tdur}} \sum_{k = -\infty}^{\infty} \int_{-1}^{0} \int_{0}^{2} e^{-\frac{(x - 4k - y)^2}{4\tdur}} \, \mathrm{d}y \, \mathrm{d}x,
\end{align*}
where $\tdur$ is the duration of contact.
The computational domain is modeled and meshed with \texttt{gmsh}\cite{gmsh}. We implemented  Algorithm~\ref{alg:Full Adaptive Algorithm} in our own \texttt{C++} code, based on the finite element library MFEM \cite{mfem:2021}. The linear systems are solved by means of the direct solver MUMPS \cite{mumps}. {The enriched finite element spaces for the goal oriented adaptivity are realized by increasing the polynomial degree of the ansatz functions by one.} The set of marked triangles is refined using red-{green} refinement \cite{banksherman1981}. Unless stated otherwise, we use a marking threshold of $0.5$. We will express convergence rates by the convention $ \mathcal{O}(h^{\alpha}) = (N_{dofs})^{-\alpha / d}$. 

Since there are no source terms present, we do not expect too much refinements for $t < \tcont $. 
At time $ t = \tcont $, there is a singularity at the spatial origin $ x = 0$, which needs to be properly resolved. This singularity is {equivalent} to the singularity of a related initial boundary value problem which starts at $ t = \tcont $ and has discontinuous initial data, with the point of discontinuity at $ x = 0$. Hence we expect heavy refinements toward $ (0, \tcont) $. 
Moreover, since our domain of interest is only on the ``right'' domain for  $ t > \tsep$, we do not expect many refinements on the ``left'' domain for $ t > \tsep$.

\begin{figure}
    \centering%
    \includegraphics[width=.95\linewidth]{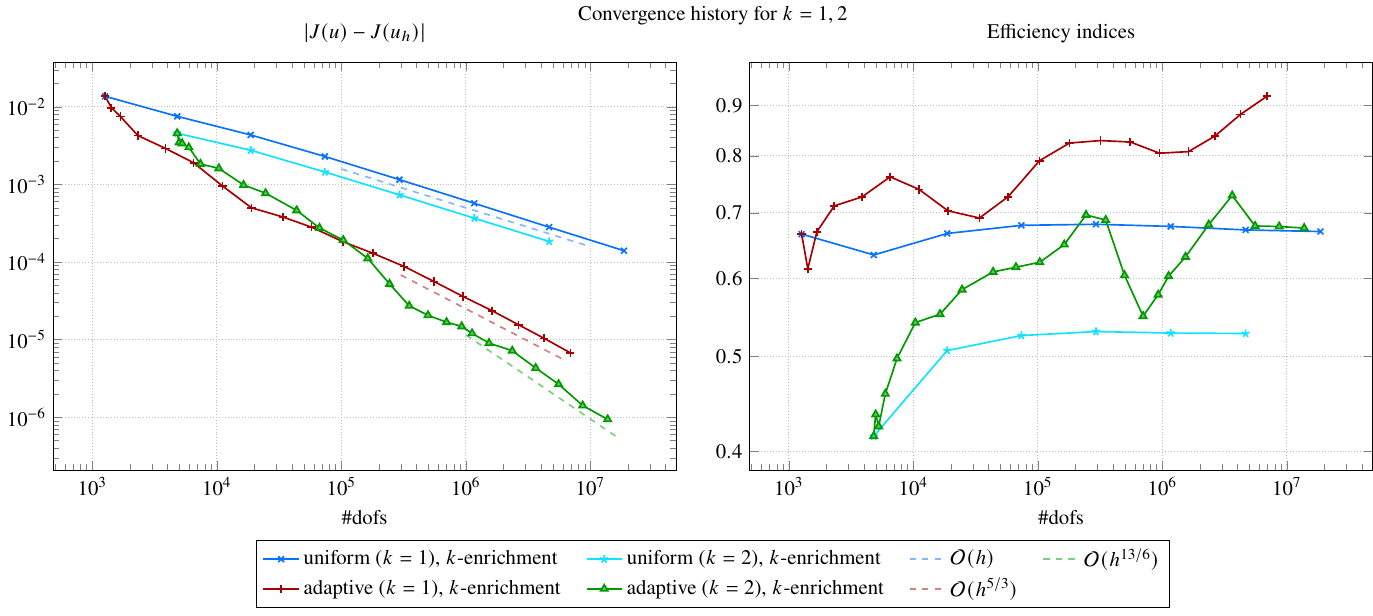}
    \caption{Convergence rates of the error in the functional and the corresponding efficiency index for the space-time domain with pre-contact part.}\label{fig:1}
\end{figure}

In Figure~\ref{fig:1}, we present on the left side the convergence history of the error in the functional, i.e.\ $ |J(u) - J(u_h)| $, and on the right the efficiency index $ \mathrm{I_{eff}} = \frac{|\eta_h|}{|J(u) - J(u_h)|}$, for different polynomial degrees $ k $. First, we note that uniform refinements always result in a rate of $ \mathcal{O}(h)$, regardless of the polynomial degree. The polynomial degree only affects the constant. Applying Algorithm~\ref{alg:Full Adaptive Algorithm} with linear ansatz functions, i.e.\ $k=1$, we observe a convergence rate of $ \mathcal{O}(h^{5/3})$, while simultaneously reducing the number of dofs needed to reach a certain error threshold by almost two orders of magnitude. In terms of efficiency, we observe that the efficiency index tends towards $ \sim0.9 $.
Increasing the polynomial degree to $ k = 2$ does improve the convergence behavior of the error in the goal functional to an observed rate of $\mathcal{O}(h^{13/6})$, but does not improve the efficiency index. 

\begin{figure}
    \centering%
    \includegraphics[width=.475\linewidth]{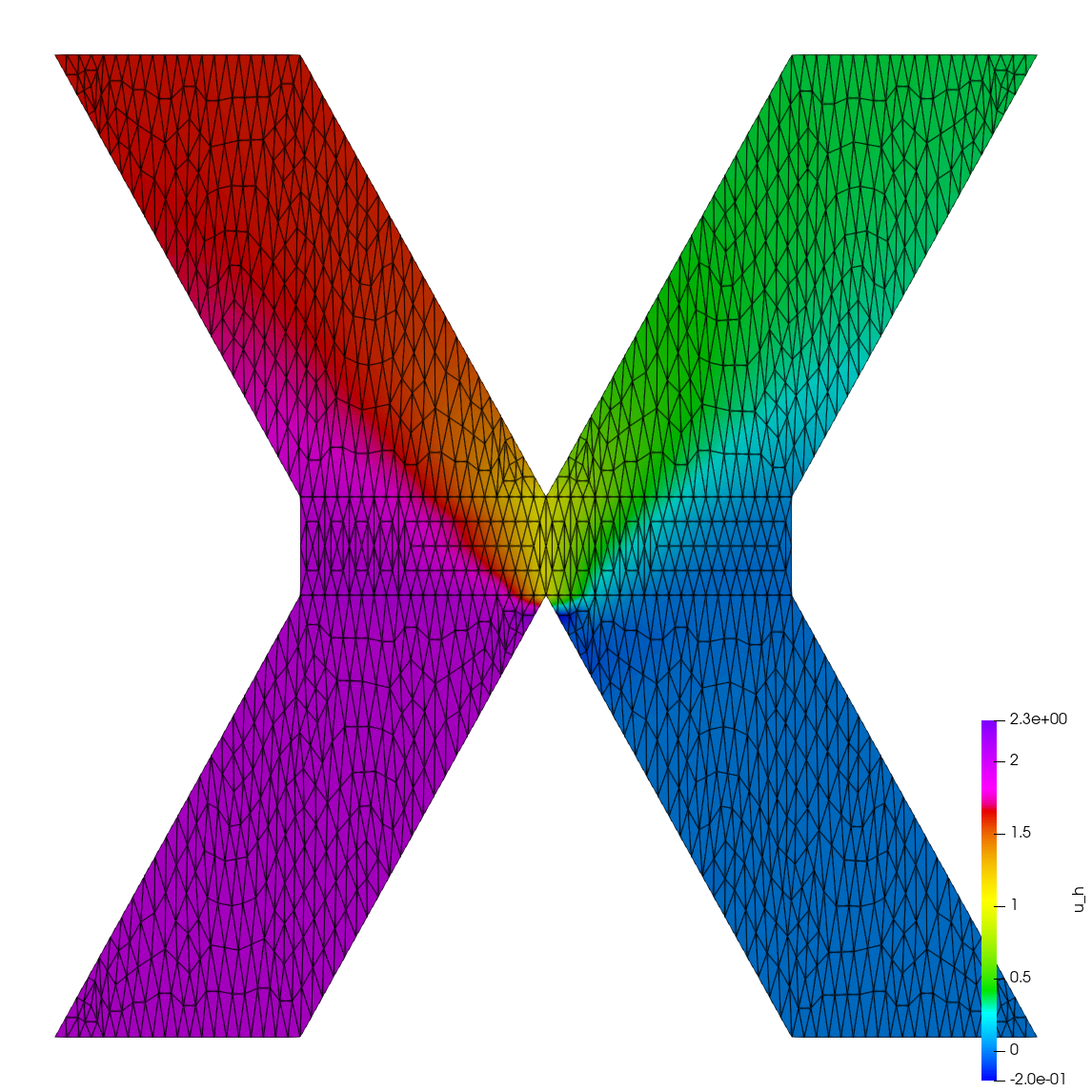}%
    \hfill%
    \includegraphics[width=.475\linewidth]{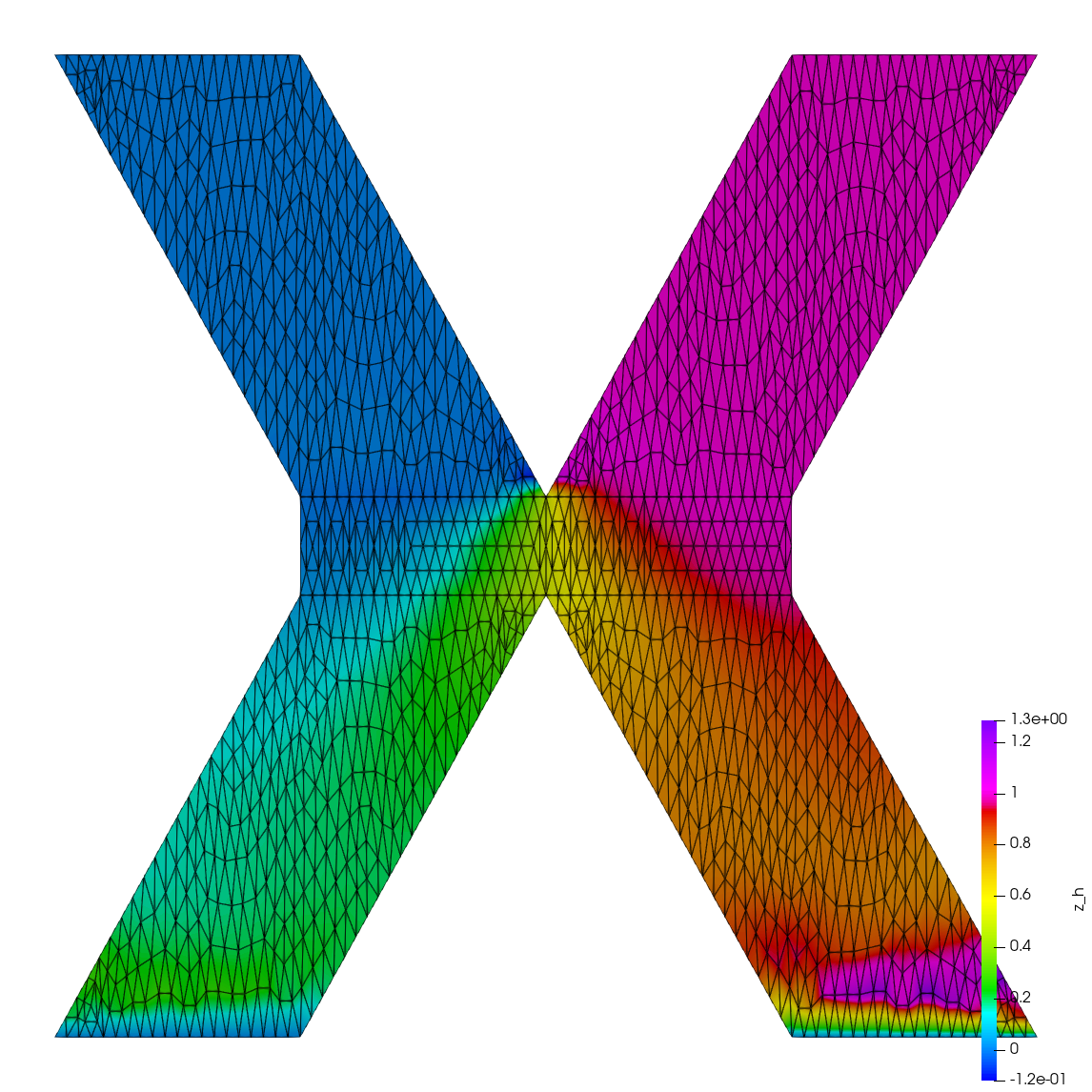}%
    \hfill%
    \includegraphics[width=.475\linewidth]{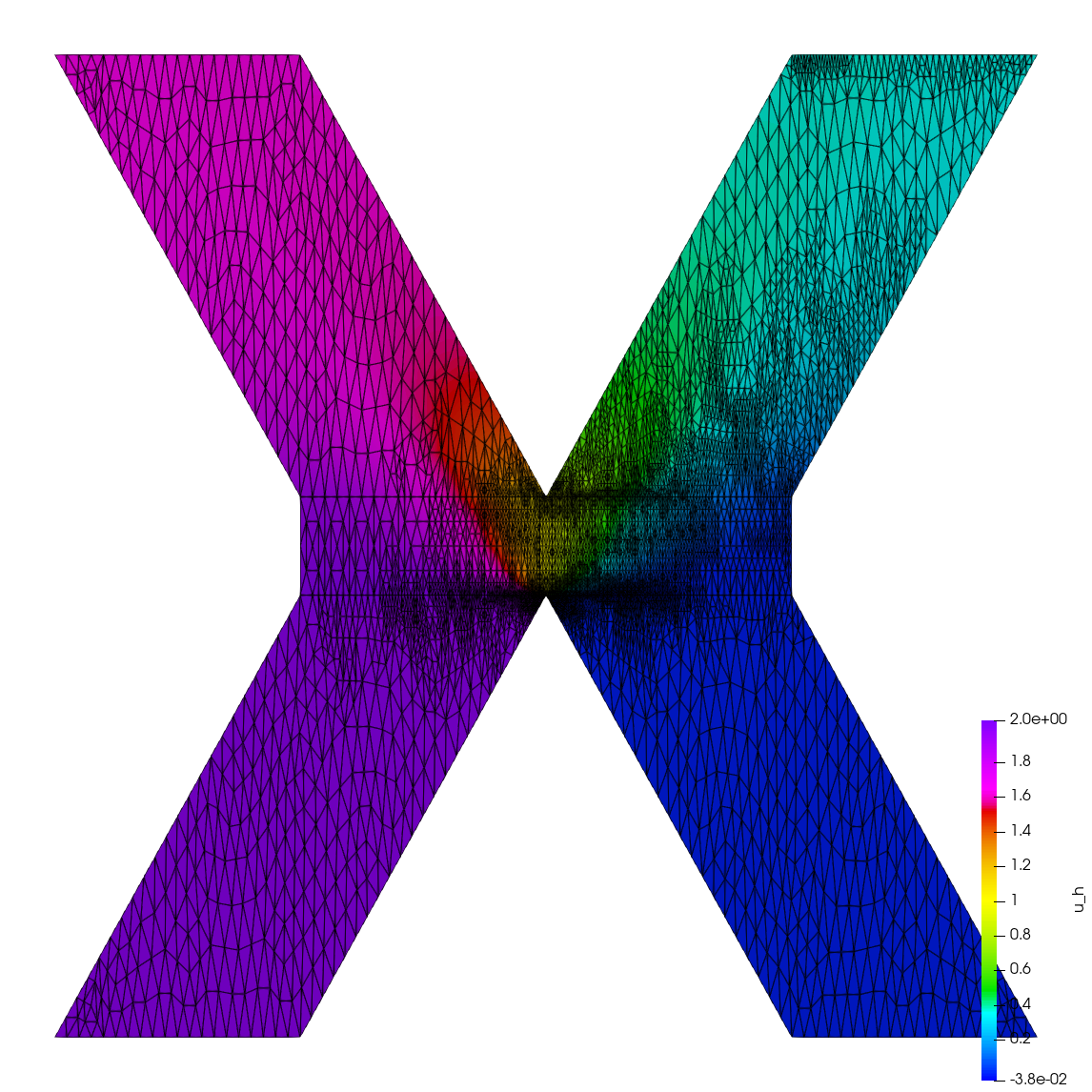}%
    \hfill%
    \includegraphics[width=.475\linewidth]{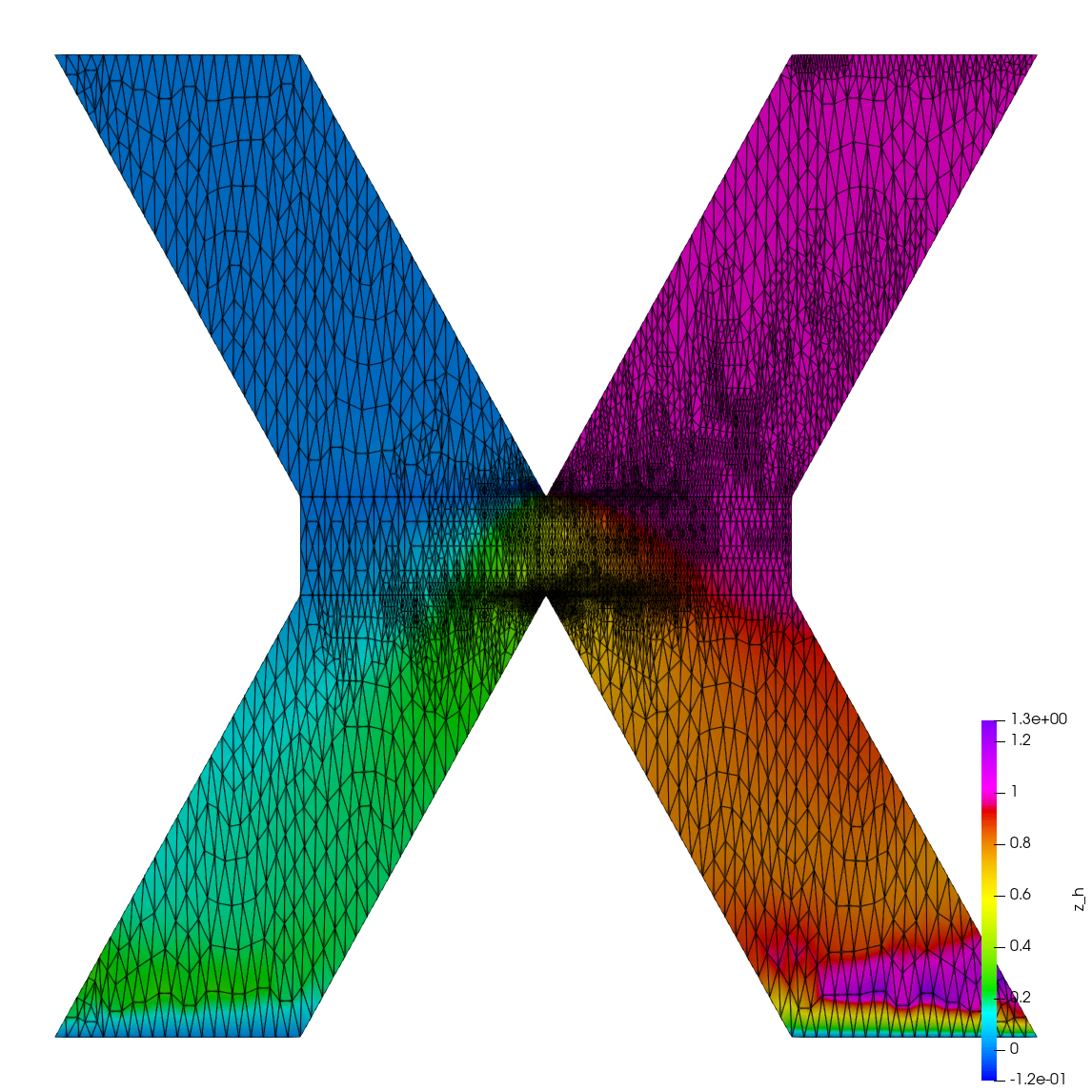}%
    \caption{Initial space-time solution $u_h$ with mesh (upper left); initial space-time sensitivity $z_h$ with mesh (upper right); space-time solution $u_h$ with mesh after {5} adaptive refinements (lower left); space-time sensitivity $z_h$ with mesh after {5} adaptive refinements (lower right); using linear finite elements and $k$-enrichment for the DWR-estimator. The mesh was stretched in $t$-direction for visualization purposes only.}%
    \label{fig:meshes}
\end{figure}

In Figure~\ref{fig:meshes}, we present the finite element solution $u_h$ and the adjoint solution (sensitivity) $ z_h$, plotted over the mesh. 
We observe that the refinements concentrate on three main regions: the region around the contact singularity $ (0,\tcont)$, the region around the separation point $ (0,\tsep)$, and around our region of interest $ \{ (x, T) \vcentcolon 1 < x < 2 \}$. 
We note that the refinements for $ t < \tcont $ on the right domain are on the one hand due to the closure step of the mesh refinement algorithm to preserve the quasi-uniformity, and on the other hand due to the loss of causality by using continuous-in-time ansatz functions. 
{While the point of separation $ (0,\tsep) $ is not a singularity for the primal problem, it {acts as} a singularity for the adjoint solution. The adjoint problem is {formally} a well posed backward heat equation, hence the point of separation {appears as} a contact point.}
The refinements towards the region of interest are as expected.

\begin{acknowledgement}
This work has been supported by
the Cluster of Excellence PhoenixD (EXC 2122, Project ID 390833453).
Furthermore, the first author is funded by an Humboldt Postdoctoral
Fellowship. Additionally, we thank Prof. Johannes Lankheit for helpful discussions.
\end{acknowledgement}

\bibliographystyle{spmpsci}
\bibliography{arxiv}

\end{document}